\documentclass[envcountsect,envcountsame]{svmult}
\usepackage{mathptmx}      
\usepackage{helvet}        
\usepackage{courier}       
\usepackage{type1cm}       
\usepackage{makeidx}        
\usepackage{graphicx}       
\usepackage{multicol}       
\usepackage[bottom]{footmisc}
\usepackage{graphicx}
\usepackage{amsmath}
\usepackage{amssymb}
\usepackage{color}
\usepackage{subcaption}
\usepackage{tikz,pgfplots}
\usepackage{wasysym}
\usepackage{wrapfig}
\usepackage{prettyref}
\usepackage{mathrsfs}
\usepackage[utf8]{inputenc}
\usepackage[
	pdftitle={An SVD in Spherical Surface Wave Tomography},
	pdfauthor={Ralf Hielscher, Daniel Potts and Michael Quellmalz},
	pdfproducer={},
	pdfcreator={}
	]{hyperref}
\newcommand{\N}{\ensuremath{\mathbb{N}_0}}
\renewcommand{\S}{\ensuremath{\mathbb{S}}}
\newcommand{\Z}{\ensuremath{\mathbb{Z}}}
\newcommand{\R}{\ensuremath{\mathbb{R}}}
\newcommand{\C}{\ensuremath{\mathbb{C}}}
\newcommand{\SO}{\ensuremath{{\mathrm{SO}(3)}}}
\newcommand{\abs}[1]{\ensuremath{\left\vert#1\right\vert}}
\newcommand{\dx}{\mathrm{d}}
\newcommand{\e}{\mathrm{e}}

\newcommand{\zb}[1]{\ensuremath{\pmb{#1}}}
\renewcommand{\d}{\, \mathrm{d}}
\renewcommand{\Box}{\hspace*{0ex} \hfill \rule{1.5ex}{1.5ex} \\}

\newcommand{\sphere}{\mathbb S}
\newcommand{\norm}[1]{\left\lVert {#1} \right\rVert}

\newcommand{\bend}{\hspace*{0ex} \hfill \hbox{\vrule height
		1.5ex\vbox{\hrule width 1.4ex \vskip 1.4ex\hrule  width 1.4ex}\vrule
		height 1.5ex}}
\newenvironment{Theorem}{\begin{theorem}}{\end{theorem}}
\newenvironment{Proof}{\goodbreak \begin{proof}}{\hspace*{0ex}\hfill\Box\end{proof}}
\newenvironment{Remark}{\begin{remark}}{\hspace*{0ex}\hfill\bend\end{remark}}
\newenvironment{Lemma}{\begin{lemma}}{\end{lemma}}

\newenvironment{Proposition}{\begin{proposition}}{\end{proposition}}
\newrefformat{def}{Definition \ref{#1}}
\newrefformat{sub}{Subsection \ref{#1}}
\newrefformat{prop}{Proposition \ref{#1}}
\newrefformat{fig}{Fig.\ \ref{#1}}
\newrefformat{rem}{Remark \ref{#1}}
\newrefformat{cor}{Corollary \ref{#1}}
\makeatletter
\long\def\theequation{\thesection.\@arabic{equation}}
\begin{document}
\title*{An SVD in Spherical Surface Wave Tomography}
\author{Ralf Hielscher, Daniel Potts and Michael Quellmalz}
\institute{Ralf Hielscher \at Chemnitz University of Technology, Faculty of Mathematics, Reichenhainer Straße 39, 09126 Chemnitz, Germany, \email{ralf.hielscher@mathematik.tu-chemnitz.de}
	\and Daniel Potts \at Same department, \email{daniel.potts@mathematik.tu-chemnitz.de}
	\and Michael Quellmalz \at Same department, \email{michael.quellmalz@mathematik.tu-chemnitz.de}
}
\date{\today}
\maketitle
\abstract{
In spherical surface wave tomography, one measures the integrals of a function defined on the sphere along great circle arcs.
This forms a generalization of the Funk--Radon transform, which assigns to a function its integrals along full great circles.
We show a singular value decomposition (SVD) for the surface wave tomography provided we have full data.\\
Since the inversion problem is overdetermined, we consider some special cases in which we only know the integrals along certain arcs.
For the case of great circle arcs with fixed opening angle, we also obtain an SVD that implies the injectivity,
generalizing a previous result for half circles in [Groemer, On a spherical integral transform and sections of star bodies, Monatsh.\ Math., 126(2):117--124, 1998].
Furthermore, we derive a numerical algorithm based on the SVD and illustrate its merchantability by numerical tests.
}

\section{Introduction}

While the famous two-dimensional Radon transform assigns to a function
$f \colon \R^{2} \to \R$ all its line integrals, its spherical generalization,
the {Funk--Radon transform} $\mathcal F \colon C(\S^2) \to C(\S^2)$, assigns
to a function on the two-dimensional sphere
$\sphere^{2} = \{\zb \xi\in \R^{3} : \abs{\zb \xi}=1\}$ its integrals
\begin{equation*}
\mathcal Ff (\zb \xi)
= \frac{1}{2\pi} \int_{\left<\zb\xi,\zb\eta\right>=0} f(\zb\eta) \d \zb\eta, \qquad
\zb{\xi}\in\S^{2},
\end{equation*}
along all great circles
$\{\zb \eta \in \sphere^{2}: \zb \eta \perp \zb \xi \}$,
$\zb \xi \in \sphere^{2}$. The investigation of the {Funk--Radon transform} dates
back to the work of Funk \cite{Fun13}, who showed the injectivity of the operator $\mathcal F$
for even functions. In other publications, the operator $\mathcal F$ is also
known as {Funk transform}, {Minkowski--Funk transform} or {spherical Radon
  transform}.

Similar to the Radon transform, the {Funk--Radon transform} plays an
important role in imaging. Motivated by
specific imaging modalities, the {Funk--Radon transform} has been generalized
further to other paths of integration, namely circles with fixed diameter
\cite{Sch69,Rub00},
circles containing the north pole \cite{AbDa93,Dah01,Hel11,Rub15II},
circles perpendicular to the equator \cite{GiReSh94,ZaSc10,HiQu15circav},
and nongeodesic hyperplane sections of the sphere
\cite{Sal15,Qu17,Pal16,Pal17}.
The integrals along half great circles have been investigated in
\cite{Gro98,GoWe06,Rub16}. Interestingly, some of these generalizations lead to
injective operators.

In this paper, we replace the great circles as paths of integration in the Funk--Radon transform by great circle arcs with arbitrary opening angle. Let
$\zb \xi,\, \zb \zeta \in \sphere^{2}$ be two points on the sphere that are
not antipodal and denote by $\gamma(\zb\xi, \zb\zeta)$ the shortest geodesic
connecting both points. Then we aim at recovering
$f \colon \sphere^{2} \to \C$ from the integrals
\begin{equation}
  \label{eq:S2xS2Radon}
  g(\zb\xi, \zb\zeta) = \int_{\gamma(\zb\xi, \zb\zeta)} f(\zb \eta) \d\zb\eta,
  \quad \zb \xi, \zb \zeta \in \sphere^{2}, \, \zb \xi \ne -\zb \zeta.
\end{equation}

The study of this problem is motivated by spherical surface wave tomography.
There, one measures the time a seismic wave travels along the Earth's surface from an
epicenter to a receiver. Knowing the traveltimes of such waves between many pairs of epicenters and receivers, one wants to recover the local phase velocity.  A common
approach is the great circle ray approximation, where it is assumed that a
wave travels along the arc of the great circle connecting epicenter and receiver. Then the traveltime of the wave equals the integral of the "slowness function"
 along the great circle arc connecting the epicenter and
the receiver, where the slowness function is defined as one over the local phase velocity \cite{WoDz84,TrWo95,No08}.  Hence, recovering the local phase
velocity as a real-valued spherical function from its mean values along
certain arcs of great circles is modeled by \eqref{eq:S2xS2Radon}, see
\cite{AmMiSi08}.

Although \eqref{eq:S2xS2Radon} uses a very intuitive parametrization of great
circle arcs on the sphere, it is not well suited for analyzing the underlying
operator since the arc length is restricted to $[0,\pi)$ and, even for
continuous $f$, the function $g$ has no continuous extension to a function on
$\sphere^{2} \times \sphere^{2}$. Therefore, we parameterize the
manifold of all great circle arcs by the arclength $2\psi \in [0,2\pi]$ and the
rotation $Q \in \mathrm{SO}(3)$ that maps the arc of integration to the
equator such that the midpoint of the arc is mapped to $(1,0,0)^{\top}$. Using
this parametrization, the arc transform is defined in Section
\ref{sec:arc-transform} as an operator
\begin{equation*}
  \mathcal A \colon C(\sphere^{2}) \to C(\SO \times [0,\pi]).
\end{equation*}
In Theorem \ref{thm:svd}, we derive a singular value decomposition of
$\mathcal A$, which involves spherical harmonics in $L^{2}(\sphere^{2})$ and
Wigner-D functions in $L^{2}(\SO \times [0,\pi])$. Furthermore, we give upper and
lower bounds for the singular values.

Since the function $f$ lives on a two-dimensional manifold but the
transformed function $\mathcal A f$ lives on a four-dimensional manifold, the
inverse problem is highly overdetermined. For this reason, we consider in
Section \ref{sec:special-cases} specific subsets of arcs that still allow for
the reconstruction of the function $f$. Most notably, we investigate in
Section \ref{sec:fixed-length} the restriction of the arc transform to arcs of
constant opening angle. This restriction includes as special cases the
ordinary Funk--Radon transform as well as the half circle transform \cite{Gro98}. For the
restricted operator, we prove in Theorem \ref{thm:mu-lim} a singular value
decomposition and show that the singular values decay as
${(n+\frac12)^{-\frac12}C(\psi,n)}$. While for opening angles $2\psi<\pi$ the
constant $C(\psi,n)$ is independent of $n$, it converges to zero for odd $n$
and $2\psi \to 2 \pi$.

Finally, we present in Section \ref{sec:numeric} a numerical algorithm for the
arc transform with fixed opening angle, which is based on the
nonequispaced fast spherical Fourier transform \cite{KeKuPo06b} and the nonequispaced fast
$\SO$ Fourier transform \cite{poprvo07}.

\section{Fourier analysis on $\S^2$ and $\SO$}
\label{sec:spherical-harmonics}

In this section, we present some basic facts about harmonic analysis on the sphere $\S^2$ and the rotation group $\SO$ and introduce the notation we will use later on.

\subsection{Harmonic analysis on the sphere}

In this section, we are going to summarize some basic facts about harmonic analysis
on the sphere as it can be found, e.g., in \cite{Mic13,DaXu2013,frgesc}.
We denote by $\Z$ the set of integers and with $\N$ the nonnegative integers.

We define the
two-dimensional sphere $\sphere^{2}=\{\zb\xi\in\R^{3}:\abs{\zb\xi}=1\}$ as
the set of unit vectors $\zb\xi = (\xi_1,\xi_2,\xi_3)^\top$ in the three-dimensional Euclidean space and make use
of its parametrization in terms of the spherical coordinates
\begin{equation*}
\zb\xi(\varphi,\vartheta)
=(\cos\varphi\, \sin\vartheta ,\sin\varphi\, \sin\vartheta ,\cos\vartheta)^{\top}
, \quad\varphi\in[0,2\pi),\; \vartheta\in[0,\pi].
\end{equation*}
Let  $f\colon\sphere^{2}\to\mathbb{C}$ be some measurable function.
With respect to spherical coordinates, the surface measure $\d{\zb\xi}$ on the sphere
reads as
\begin{equation*}
\int_{\sphere^{2}}f(\zb\xi)\d{\zb\xi}
=\int_{0}^{\pi}\int_{0}^{2\pi}f(\zb\xi(\varphi,\vartheta))\, \sin\vartheta \d{\varphi} \d{\vartheta}.
\end{equation*}
The Hilbert space $L^{2}(\sphere^{2})$ is the space of all measurable functions
$f\colon\sphere^{2}\to\mathbb{C}$ whose norm
$\norm{f}_{L^{2}(\S^2)}=({\left<f,f\right>_{L^2(\S^2)}})^{1/2}$ is finite,
where $\left<f,g\right>_{L^2(\S^2)} = \int_{\sphere^{2}} f(\zb\xi)\,\overline{g(\zb\xi)} \d\zb\xi$ denotes the usual $L^2${\textendash}inner product.

We define the associated Legendre functions
\begin{equation*}
P_{n}^{k}(t)=\frac{(-1)^{k}}{2^{n}n!}\left(1-t^{2}\right)^{k/2}\frac{\mathrm{d}^{n+k}}{\mathrm{d}t^{n+k}}\left(t^{2}-1\right)^{n}
,\quad t\in[-1,1],
\end{equation*}
of degree $n\in\N$ and order $k=0,\dots,n$.
We define the normalized associated Legendre functions by
\begin{equation}
\widetilde P_n^k
= \sqrt{\frac{2n+1}{4\pi}\,\frac{(n-k)!}{(n+k)!}}\, P_{n}^{k}
\label{eq:associated-legendre-normalized}
\end{equation}
and
$$\widetilde P_n^{-k} = (-1)^k\, \widetilde P_n^k,$$
where the factor $(-1)^k$ is called Condon--Shortley phase, which is omitted by some authors.

An orthonormal basis in the Hilbert space $L^{2}(\sphere^{2})$ of square integrable
functions on the sphere is formed by the spherical harmonics
\begin{equation}
Y_{n}^{k}(\zb\xi(\varphi,\vartheta))=
\widetilde P_n^k(\cos\vartheta)\, \mathrm{e}^{\mathrm{i}k\varphi}
\label{eq:spherical-harmonics}
\end{equation}
of degree $n\in\N$ and order $k=-n,\dots,n$.
Accordingly, any function $f\in L^{2}(\sphere^{2})$ can be expressed by its spherical Fourier series
\[
 f= \sum_{n=0}^{\infty} \sum_{k=-n}^{n} \hat{f}_n^k\, Y_{n}^{k}
 \]
 with the spherical Fourier coefficients
\[
\hat{f}_n^k=\int_{\sphere^{2}}f(\zb\xi)\,\overline{Y_{n}^{k}(\zb\xi)}\d{\zb\xi},\quad n\in\N,\;k=-n,\dots,n.
\]
We define the space of spherical polynomials of degree up to $N\in\N$ by
\begin{equation*}
\mathscr P_N
= \operatorname{span} \left\{ Y_n^k : n=0,\dots,N,\, k=-n,\dots,n \right\}.
\end{equation*}

\subsection{Rotational harmonics}
We state some facts about functions on the rotation group $\SO$.
This introduction is based on \cite{HeHeKi98}, rotational Fourier transforms date back to Wigner, 1931, see \cite{Wig31}.
The rotation group $\SO$ consists of all orthogonal $3\times 3$--matrices with determinant one equipped with the matrix multiplication as group operation.
Every rotation $Q\in \SO$ can be expressed in terms of its Euler angles $\alpha,\beta,\gamma$ by
\begin{equation*}
Q(\alpha,\beta,\gamma)
= R_3 (\alpha) R_2 (\beta) R_3 (\gamma),
\quad \alpha, \gamma \in [0,2\pi),\; \beta \in [0,\pi],
\end{equation*}
where $R_i (\alpha)$ denotes the rotation of the angle $\alpha$ about the $\zb{\xi}_i$--axis, i.e.,
\begin{align*}
R_3(\alpha) &=
\begin{pmatrix}
\cos\alpha &-\sin\alpha &0\\
\sin\alpha &\cos\alpha &0\\
0&0&1
\end{pmatrix},
&
R_2(\beta) &=
\begin{pmatrix}
\cos\beta&0 &\sin\beta \\
0&1&0\\
-\sin\beta&0 &\cos\beta
\end{pmatrix}.
\end{align*}
Note that we use this $zyz$-convention of Euler angles throughout this paper.
The integral of a function $g \colon \SO\to\C$ on the rotation group is given by
\begin{equation*}
\int_{\SO} g(Q) \d Q
= \int_{0}^{2\pi} \int_{0}^{\pi} \int_{0}^{2\pi}  g(Q(\alpha,\beta,\gamma)) \sin (\beta) \d \alpha \d \beta \d \gamma.
\end{equation*}
We define the rotational harmonics or Wigner D-functions
$D_n^{k,j}$ of degree $n\in\N$ and orders $k,j \in \{-n,\dots,n\}$ by
\begin{equation*}
D_n^{k,j} (Q(\alpha,\beta,\gamma))
= \e^{-\mathrm{i}k\alpha} d_n^{k,j} (\cos\beta) \e^{-\mathrm{i}j\gamma},
\end{equation*}
where the Wigner d-functions are given by \cite[p.\ 77]{Varsha88}
\begin{align*}
 d_n^{k,j}(t)
&= \frac{(-1)^{n-j}}{2^n} \sqrt\frac{(n+k)!(1-t)^{j-k}}{(n-j)!(n+j)!(n-k)!(1+t)^{j+k}}
\left(\frac{\dx}{\dx t}\right)^{n-k}\frac{(1+t)^{n+j}}{(1-t)^{-n+j}}
.
\end{align*}
The Wigner d-functions satisfy the orthogonality relation
\begin{equation*}
\int_{-1}^{1} d_n^{k,j} (t)\, d_{n'}^{k,j} (t) \d t
= \frac{2\delta_{n,n'}}{2n+1}.
\end{equation*}
We define the space of square-integrable functions $L^2(\SO)$ with inner product $\left<f,g\right>_{L^2(\SO)}=\int_\SO f(Q)\,\overline{g(Q)}\d Q$.
By the Peter--Weyl theorem, the rotational harmonics $D_n^{k,j}$ are complete in $L^2(\SO)$ and satisfy the orthogonality relation
\begin{equation}
\left< D_n^{k,j} ,D_{n'}^{k',j'}\right>_{L^2(\SO)}=
\int_{\SO} D_n^{k,j}(Q) \,\overline{D_{n'}^{k',j'}(Q)} \d Q
= \frac{8\pi^2}{2n+1}\, \delta_{n,n'} \delta_{k,k'} \delta_{j,j'}.
\label{eq:rotational-harmonics-orthogonal}
\end{equation}
We define the rotational Fourier coefficients of $g\in L^2(\SO)$ by
\begin{equation}
\hat g_n^{k,j}
= \frac{2n+1}{8\pi^2} \left<g,D_n^{k,j}\right>_{L^2(\SO)}
,\quad n\in\N,\; k,j = -n,\dots,n.
\label{eq:Fourier-coefficients-SO3}
\end{equation}
Then the rotational Fourier expansion of $g$ holds
\begin{equation*}
g = \sum_{n=0}^{\infty} \sum_{k,j=-n}^{n} \hat g_n^{k,j}\, D_n^{k,j}.
\end{equation*}
The rotational Fourier transform is also known as $\SO$ Fourier transform (SOFT) or Wigner D-transform.

The rotational harmonics $D_n^{k,j}$ are eigenfunctions of the Laplace--Beltrami operator on $\SO$ with the corresponding eigenvalues $-n(n+1)$.
The rotational harmonics $D_n^{j,k}$ are the matrix entries of the left regular representations of $\SO$, see \cite{hiediss07,Varsha88}.
In particular, the rotation of a spherical harmonic satisfies
\begin{equation}
\label{eq:spherical-harmonic-rotation}
Y_n^k(Q^{-1}\zb\xi)
= \sum_{j=-n}^{n} D_n^{j,k}(Q)\, Y_n^j({\zb\xi}).
\end{equation}

\subsection{Singular value decomposition}

Let $\mathcal K\colon X\to Y$ be a compact linear operator between the separable Hilbert spaces $X$ and $Y$.
A singular system
$ \{(u_n,v_n,\sigma_n): n\in\N\}$
consists of an orthonormal basis $\{u_n\}_{n=0}^\infty$ of $X$, an orthonormal basis $\{v_n\}_{n=0}^\infty$ in the closed range of $\mathcal K$ and singular values $\sigma_n\to0$
such that operator $\mathcal K$ can be diagonalized as
\begin{equation*}
\mathcal Kx = \sum_{n=0}^{\infty} \sigma_n \left<x,u_n\right> v_n
,\qquad x\in X.
\end{equation*}
If all singular values $\sigma_n$ are nonzero, the operator $\mathcal K$ is injective and for $y=\mathcal Kx$, we have
$$
x = \sum_{n=0}^{\infty} \frac{\left<y,v_n\right>}{\sigma_n} u_n.
$$
The instability of an inverse problem can be characterized by the decay of the singular values.
The problem of solving $\mathcal Kx=y$ for $x$ is called mildly ill-posed of degree $\alpha>0$ if $\sigma_n \in \mathcal O(n^{-\alpha})$, cf.\ \cite[Sec.\ 2.2]{EnHaNe96}.

\begin{wrapfigure}[11]{r}{0.5\textwidth}
	\vspace{-5em}
	\includegraphics[width=0.49\textwidth]{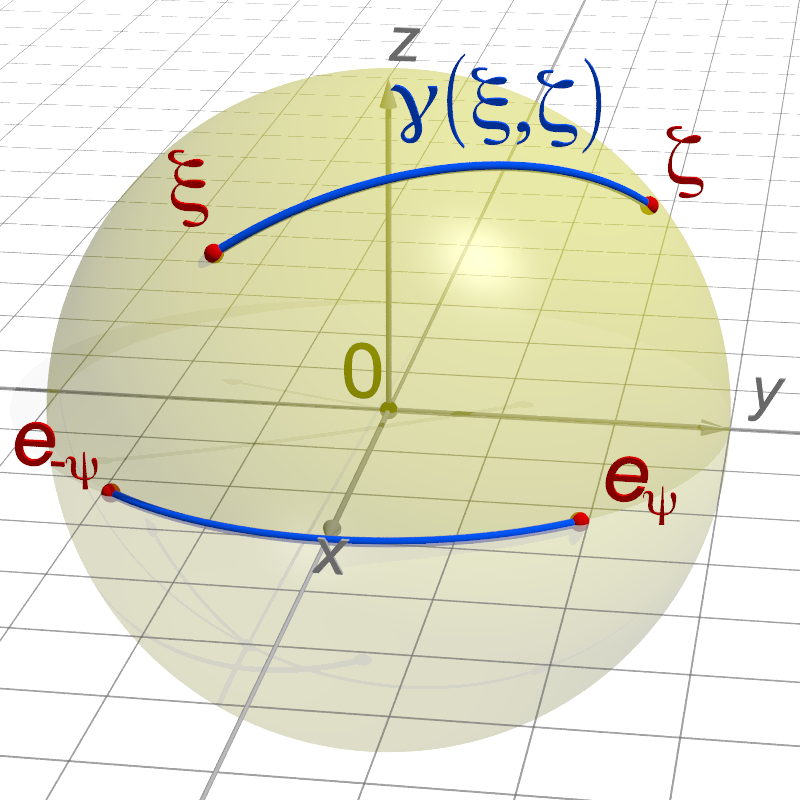}
	\caption{Visualiation of the arc $\gamma(\zb\xi,\zb\zeta)$}
	\label{fig:circle-arc}
\end{wrapfigure}

\section{Circle arcs}

For any two points $\zb\xi, \zb\zeta$  on the sphere $\S^2$ that are not antipodal, there exists a shortest geodesic $\gamma(\zb\xi, \zb\zeta)$ between $\zb\xi$ and $\zb\zeta$.
This geodesic is an arc of the great circle that contains $\zb\xi$ and $\zb\zeta$.

The manifold of all great circle arcs is four-dimensional since they are determined by two points $\zb\xi,\zb\zeta\in\S^2$ and only coincide when $\zb\xi$ and $\zb\zeta$ are interchanged.

\subsection{The arc transform}
\label{sec:arc-transform}

A great circle arc $\gamma(\zb\xi,\zb\zeta)$ can be parameterized by its length
$2 \psi
= \arccos(\left<\zb\xi,\zb\zeta\right>)$ and a
rotation $Q \in \SO$ which is defined as follows.
Let
$$
\zb e_\varphi = (\cos\varphi,\,\sin\varphi,\,0)^\top \in \S^2
$$
be the point on the equator of $\S^2$ with latitude $\varphi\in\R$.
Then there exists a unique rotation $Q\in \SO$ such that $Q(\zb\xi) = \zb e_{-\psi}$ and $Q(\zb\zeta) = \zb e_{\psi}$.
Such an arc $\gamma$ and its rotation are depicted in \prettyref{fig:circle-arc}.
With this definition, the integral over the arc $\gamma(\zb\xi,\zb\zeta)$ may
be rewritten as
\begin{equation*}
\int_{\gamma(\zb\xi,\zb\zeta)} f(\zb\eta) \d\zb\eta
= \int_{Q\gamma(\zb\xi,\zb\zeta)} f(Q^{-1}\zb\eta) \d\zb\eta
= \int_{-\psi}^{\psi} f\circ Q^{-1} (\zb e_{\varphi}) \d\varphi.
\end{equation*}
This motivates the following definition of the arc transform
\begin{equation}
\begin{split}
\mathcal A\colon C(\S^2)&\to C(\SO\times[0,\pi]),
\\
\mathcal Af(Q,\psi)
&= \int_{-\psi}^{\psi} f\circ Q^{-1} (\zb e_{\varphi}) \d\varphi.
\end{split}
\label{eq:a}
\end{equation}
The great circle arcs $\gamma(\zb\xi,\zb\zeta)$ and $\gamma(\zb\zeta,\zb\xi)$ are identical.
This symmetry also holds for the operator $\mathcal A$.
Using Euler angles, we have the identity
\begin{equation*}
\mathcal A f(Q(\alpha,\beta,\gamma),\psi)
= \mathcal A f(Q(2\pi-\alpha,\pi-\beta,\gamma+\pi),\psi),
\end{equation*}
where we assume the Euler angle $\gamma$ as $2\pi$-periodic.

\subsection{Singular value decomposition of the arc transform}

In the following, we use double factorials defined by
$n!!=n(n-2)\cdots 2$ for $n$ even or $n!!=n(n-2)\cdots 1$ for $n$ odd and $0!! = (-1)!!=1$.
The next theorem shows how the arc transform $\mathcal A$ acts on spherical harmonics $Y_n^k$.
The corresponding result for the parametrization in terms of the endpoints of an arc is found in \cite[Appx.\ C]{DaTr98}, see also \cite{AmMiSi08}.

\begin{Theorem}
	\label{thm:arc-transform-sh}
	Let $n\in\N$ and $k\in\{-n,\dots,n\}$.
	Then
	\begin{equation}
	\mathcal A Y_n^k (Q,\psi)
	= \sum_{j=-n}^n \widetilde P_{n}^{j}(0)\, D_n^{j,k}(Q)\, s_j(\psi),
	\label{eq:AY}
	\end{equation}
	where
	\begin{equation}
	s_j(\psi)=
	\begin{cases}
	2\psi, &j=0\\
	\frac{2\sin(j\psi)}{j}, &j\neq0
	\end{cases}
	\label{eq:sj}
	\end{equation}
	and
	\begin{equation}
	\label{eq:Pnj0}
	\widetilde P_n^j(0)
	= \begin{cases}
	(-1)^{\frac{n+j}{2}} \sqrt{\frac{2n+1}{4\pi} \frac{(n-j-1)!!(n+j-1)!!}{(n-j)!!(n+j)!!}}, & n+j\text{ even}\\
	0, &n+j\text{ odd}.
	\end{cases}
	\end{equation}
\end{Theorem}
\begin{Proof}
By \eqref{eq:spherical-harmonic-rotation}, we obtain
\begin{align*}
\mathcal A Y_n^k (Q,\psi)
= \int_{-\psi}^{\psi} Y_n^k (Q^{-1} (\zb e_{\varphi}))
= \sum_{j=-n}^n D_n^{j,k}(Q) \int_{-\psi}^\psi Y_n^j (\zb e_{\varphi}) \d\varphi.
\end{align*}
By the definition \eqref{eq:spherical-harmonics} of the spherical harmonics, we see that
\begin{align*}
\int_{-\psi}^{\psi} Y_n^j (\zb e_{\varphi}) \d\varphi
=\widetilde P_{n}^{j}(0) \int_{-\psi}^{\psi} \mathrm{e}^{\mathrm{i}j\varphi} \d\varphi
=\widetilde P_{n}^{j}(0)\, s_j(\psi).
\end{align*}
Hence,
\begin{equation*}
\mathcal A Y_n^k (Q,\psi)
= \sum_{j=-n}^n D_n^{j,k}(Q)\, \widetilde P_{n}^{j}(0)\, s_j(\psi).
\end{equation*}
Now we calculate $\widetilde P_n^j(0)$.
By \cite{HiQu15circav}, $P_n^j(0) = 0$ if $n+j$ is odd and otherwise
\begin{equation*}
P_n^j(0)
= (-1)^{\frac{n+j}{2}}\frac{(n+j-1)!!}{(n-j)!!}.
\end{equation*}
Hence, we obtain by \eqref{eq:associated-legendre-normalized}
\begin{equation*}
  \begin{split}
    \widetilde P_n^j(0) &
    = \sqrt{\frac{2n+1}{4\pi}\frac{(n-j)!}{(n+j)!}}\,
    (-1)^{\frac{n+j}{2}}\frac{(n+j-1)!!}{(n-j)!!}
    \\&= (-1)^{\frac{n+j}{2}}
    \sqrt{\frac{2n+1}{4\pi} \frac{(n-j-1)!!(n+j-1)!!}{(n-j)!!(n+j)!!}}
  \end{split}
\end{equation*}
if $n+j$ is even, which implies \eqref{eq:Pnj0}.
\end{Proof}

\begin{Lemma}
	\label{lem:double-factorial}
	Let $n\in\N$ and $j\in\{-n,\dots,n\}$.
	If $n+j$ is odd, then $\widetilde P_n^j(0) = 0$.
	Otherwise, we have
	\begin{equation}
	\frac{2n+1}{2\pi^2 \sqrt{(n+1)^2-j^2}}
	\leq \abs{\widetilde P_n^j(0)}^2
	\le \frac{2n+1}{4\pi \sqrt{(n+1)^2-j^2}}.
	\label{eq:associated-legendre-bound}
	\end{equation}
	Furthermore, for $j\in\N$,
	\begin{equation}
	\lim_{\substack{n\to\infty\\n+j \text{ even}}} \abs{\widetilde P_{n}^j(0)} = \frac1\pi.
	\label{eq:Pnj0-lim}
	\end{equation}
\end{Lemma}
\begin{Proof}
	We first show that for $m\in\mathbb N$,
	\begin{equation}
	\sqrt{\frac{2}{\pi(2m+1)}} \le
	\frac{(2m-1)!!}{(2m)!!} \le \frac{1}{\sqrt{2m+1}}.
	\label{eq:wallis}
	\end{equation}
	With the definition
	\begin{equation*}
	u(m)=
	\left(\frac{(2m)!!}{(2m-1)!!}\right)^2\frac{1}{2m+1},
	\qquad m\in\N,
	\end{equation*}
	we see that $u(0)=1$ and $u$ is increasing because of $m\ge1$ and
	\begin{equation*}
	\frac{u(m)}{u(m-1)}
	= \frac{(2m)^2}{(2m-1)^2} \frac{2m-1}{2m+1}
	= \frac{(2m)^2}{(2m)^2-1}
	>1.
	\end{equation*}
	That implies the right inequality of \eqref{eq:wallis}.
	Furthermore, Wallis' product states the convergence
	\begin{equation}
	u(m) = \frac21\,\frac23\,\frac43\,\frac45\,\frac65\,\frac67\,\cdots\,\frac{2m}{2m-1}\,\frac{2m}{2m+1} \longrightarrow \frac\pi2
	\label{eq:wallis-orig}
	\end{equation}
	for $m\to\infty$, see also \cite{Bau07}.
	This shows the left inequality of \eqref{eq:wallis}.

        By \eqref{eq:Pnj0} and \eqref{eq:wallis}, we obtain the upper bound
        \begin{equation*}
          \abs{\widetilde P_n^j(0)}^2
          = \frac{2n+1}{4\pi} \frac{(n-j-1)!!}{(n-j)!!}
          \frac{(n+j-1)!!}{(n+j)!!}
          \le \frac{2n+1}{4\pi} \frac1{\sqrt{n-j+1}} \frac1{\sqrt{n+j+1}}.
        \end{equation*}
        The lower bound follows analogously.
    Moreover, we have
    \begin{align*}
    \abs{\widetilde P_{j+2m}^j(0)}^2
    &= \frac{2(j+2m)+1}{4\pi} \frac{(2m-1)!!}{(2m)!!}
    \frac{(2m+2j-1)!!}{(2m+2j)!!}
    .
    \end{align*}
    Hence, Wallis product \eqref{eq:wallis-orig} shows that for $j\in\N$
    \begin{align*}
    \lim_{m\to\infty} \abs{\widetilde P_{j+2m}^j(0)}^2
    &= \lim_{m\to\infty} \frac{2(j+2m)+1}{4\pi}
    \frac2\pi
    \frac{1}{\sqrt{2m+1}\,\sqrt{2m+2j+1}}
    \\&= \lim_{m\to\infty} \frac{2m+j+\frac12}{\pi^2\,\sqrt{(2m+j+1)^2-j^2}}
    = \frac{1}{\pi^2},
    \end{align*}
    which proves the assertion.
\end{Proof}

Next, we derive a singular value decomposition for the spherical arc
transform. To this end, we define for $n\in \N$ and $k = -n,\dots,n$ the
functions $E_n^k \in L^{2}(\SO\times[0,\pi])$ by
\begin{equation}
  \label{eq:E}
  E_n^k (Q,\psi)
  =   \sum_{j=-n}^n D_n^{j,k}(Q)\, \widetilde P_{n}^{j}(0)\, s_j(\psi),
  \qquad Q\in \SO,\;\psi\in[0,\pi].
\end{equation}

\begin{Theorem}
  \label{thm:svd}
  The operator $\mathcal A\colon L^2(\S^2)\to L^2(\SO\times[0,\pi])$ is compact with the singular value decomposition
  \begin{equation*}
    \left\{ \left(Y_n^k,\; \widetilde E_n^k,\; \sigma_n\right): n\in\N,\;k\in\{-n,\dots,n\} \right\},
  \end{equation*}
  with the singular values
  \begin{equation}
    \sigma_n = \norm{E_n^k}_{L^2(\SO\times[0,\pi])}
    =\sqrt{\frac{32\pi^3}{2n+1}} \sqrt{\frac{\pi^2}{3} \abs{\widetilde
        P_{n}^{0}(0)}^2 + \sum_{\substack{j=1}}^{n} \frac{1}{j^2}
      \abs{\widetilde P_{n}^{j}(0)}^2 }
  \label{eq:sv}
  \end{equation}
  satisfying
  \begin{align}
	\sqrt{\frac{16}{3}\pi^3}
    &\le \sigma_n \,\sqrt{n+1} \le
    \sqrt{\frac{8}{3}\pi^4 +4\pi^2} , && n\text{ even,}
    \label{eq:sv-bound-even}\displaybreak[0]\\
    4\sqrt{\pi}
    &\le \sigma_n \,\sqrt{n+1} \le
    2\pi\sqrt{\frac{4}{\sqrt{3}} + 1}, && n\text{ odd,}
    \label{eq:sv-bound-odd}
  \end{align}
  and the orthonormal function system
  $\widetilde{E}^n_k = \sigma_n^{-1} E_n^k$, $n\in\N$, $k\in\{-n,\dots,n\}$ in
  $L^2(\SO\times[0,\pi])$.
\end{Theorem}
\begin{Proof}
By the orthogonality \eqref{eq:rotational-harmonics-orthogonal} of the rotational harmonics, we have
\begin{align*}
&\hspace{-1em}\left<E_n^k, E_{n'}^{k'}\right>_{L^2(\SO\times[0,\pi])}
\\
&= \sum_{j=-n}^n \sum_{j'=-n'}^{n'} \widetilde P_{n}^{j}(0)\, \widetilde P_{n'}^{j'}(0)
\int_{\SO} D_{n}^{j,k}(Q)\, \overline{D_{n'}^{j',k'}(Q)} \d Q
\int_{0}^{\pi} s_j(\psi)\, s_{j'}(\psi) \d\psi
\\&=
\sum_{j=-n}^n \sum_{j'=-n'}^{n'} \frac{8\pi^2}{2n+1}\delta_{nn'}\,\delta_{kk'}\,\delta_{jj'}\, \widetilde P_{n}^{j}(0)\, \widetilde P_{n'}^{j'}(0) \int_{0}^{\pi} s_j(\psi)\, s_{j'}(\psi) \d\psi
\\&=
\delta_{nn'}\,\delta_{kk'} \sum_{j=-n}^n \frac{8\pi^2}{2n+1} \abs{\widetilde P_{n}^{j}(0)}^2 \int_{0}^{\pi} s_j(\psi)^2 \d\psi
\\&=
\delta_{nn'}\,\delta_{kk'} \frac{8\pi^2}{2n+1} \sum_{j=-n}^n \abs{\widetilde P_{n}^{j}(0)}^2
\begin{cases}
\frac{4\pi^3}{3}, &j=0\\
\frac{2\pi}{j^2}, &j\neq0.
\end{cases}
\end{align*}
This shows that the functions $\mathcal A Y_n^k$ are orthogonal in the space $L^2(\SO\times[0,\pi])$
and have the norm
\begin{align*}
\norm{E_n^k}^2_{L^2(\SO\times[0,\pi])}
&= \frac{8\pi^2}{2n+1} \sum_{j=-n}^n \abs{\widetilde P_{n}^{j}(0)}^2
\begin{cases}
\frac{4\pi^3}{3}, &j=0\\
\frac{2\pi}{j^2}, &j\neq0.
\end{cases}
\\
&= \frac{16\pi^3}{2n+1} \left(\frac{2\pi^2}{3} \abs{\widetilde P_{n}^{0}(0)}^2 + 2\sum_{\substack{j=1}}^{n} \frac{1}{j^2} \abs{\widetilde P_{n}^{j}(0)}^2 \right),
\end{align*}
where we used that $ \abs{\widetilde P_{n}^{j}(0)} = \abs{\widetilde P_{n}^{-j}(0)}$.
In order to prove that $\mathcal A$ is compact, we show that the singular values $\sigma_n$ decay for $n\to\infty$.
We have by Lemma \ref{lem:double-factorial} for $n=2m$ even
\begin{align*}
\sigma_{2m}^2
&\le 4\pi^2 \left(\frac{2\pi^2}{3} \frac{1}{2m+1}
+ 2\sum_{j=1}^{m} \frac{1}{(2j)^2} \frac{1}{\sqrt{(2m+1)^2-(2j)^2}}\right) .
\end{align*}
Replacing the sum by an integral, we estimate for $n$ even
\begin{align*}
2\sum_{j=1}^{m} \frac{1}{(2j)^2} \frac{1}{\sqrt{(2m+1)^2-(2j)^2}}
&\le 2\int_{1/2}^{m+1/2} \frac{1}{(2j)^2} \frac{1}{\sqrt{(2m+1)^2-(2j)^2}} \d j
\\
&= 2\left[ -\frac{\sqrt{(2m+1)^2-(2j)^2}}{2j(2m+1)^2} \right]_{1/2}^{m+1/2}
\\
&= 2\frac{\sqrt{m^2+m}}{(2m+1)^2}
\le \frac{1}{2m+1} ,
\end{align*}
where we made use of the convexity of the integrand.
Hence,
$$
\sigma_{2m}^2
\le 4\pi^2 \left(\frac{2\pi^2}{3} \frac{1}{2m+1}
+ \frac{1}{2m+1}\right)
= 4\pi^2 \left(\frac{2\pi^2}{3} +1 \right)\frac{1}{2m+1}.
$$
For odd $n=2m-1$, we proceed analogously.
We have
\begin{align*}
\sigma_{2m-1}^2
&\le 8\pi^2 \sum_{j=1}^{m} \frac{1}{(2j-1)^2} \frac{1}{\sqrt{(2m)^2-(2j-1)^2}} .
\end{align*}
Note that, for the estimation of the sum by an integral, we extract the summand for $j=1$
\begin{align*}
\sigma_{2m-1}^2
&\le 8\pi^2\left(\frac{1}{\sqrt{(2m)^2-1}} + \int_{1}^{m+1/2} \frac{1}{(2j-1)^2} \frac{1}{\sqrt{(2m)^2-(2j-1)^2}} \d j\right)
\\
&= 8\pi^2\left(\frac{1}{\sqrt{(2m)^2-1}} + \frac{\sqrt{(2m)^2-1}}{2(2m)^2} \right)
\\
&\le 8\pi^2\left(\frac{2}{\sqrt{3}\,2m} + \frac{1}{2(2m)} \right)
= 4\pi^2\left(\frac{4}{\sqrt{3}} + 1 \right)\frac{1}{2m}.
\end{align*}
For the lower bound of the singular values, we also use Lemma \ref{lem:double-factorial}.
For even $n$, we extract the summand $j=0$ and obtain
\begin{align*}
\sigma_n^2
&=\frac{16\pi^3}{2n+1} \left(\frac{2\pi^2}{3} \abs{\widetilde P_{n}^{0}(0)}^2 + 2\sum_{\substack{j=1}}^{n} \frac{1}{j^2} \abs{\widetilde P_{n}^{j}(0)}^2 \right)
\\&
\ge \frac{32\pi^5}{3(2n+1)} \abs{\widetilde P_{n}^{0}(0)}^2
\ge \frac{16\pi^3}{3(n+1)}.
\end{align*}
For odd $n$, we extract the summand $j=1$ and obtain
\begin{equation*}
\sigma_n^2
\ge \frac{32\pi^3}{2n+1}\abs{\widetilde P_{n}^{1}(0)}^2
\ge \frac{16\pi}{\sqrt{(n+1)^2-1}}
\ge \frac{16\pi}{n+1}.
\end{equation*}
\end{Proof}

The singular values $\sigma_n$ decay with rate $n^{-1/2}$.
This is the same asymptotic decay rate as of the eigenvalues of the Funk--Radon transform, cf.\ \cite{Str81}.

\section{Special cases}

\label{sec:special-cases}

The recovery of a function $f$ from the arc integrals $\mathcal Af$ is overdetermined considered we have full data.
In the following subsections, we are going to examine some special cases, where we can reconstruct $f$ from integrals only along certain arcs.

\subsection{Arcs starting in a fixed point}
\label{sec:arcs-north-pole}

As a simple example, we fix one endpoint of the arcs.
Without loss of generality, we assume that this endpoint is the north pole.
The arc connecting the north pole $\zb e^3$ and an arbitrary other point $\zb\xi(\varphi,\vartheta)\in\S^2$ is given by
$$ \gamma(\zb e^3, \zb\xi(\varphi,\vartheta))
= \{\zb\eta(\varphi,\varrho)\in\S^2: \varrho\in[0,\vartheta] \}.
$$
Since, with $Q = Q\left(\tfrac{\vartheta}{2}, \tfrac{\pi}{2},
  \tfrac{3\pi}{2}-\varphi\right)\in\SO$, we have $Q \zb e^{3} = \zb
e_{\frac{\vartheta}2}$ and $Q \zb \xi = \zb
e_{-\frac{\vartheta}2}$. The restriction  $\mathcal B
\colon C(\sphere^{2}) \to C(\sphere^{2})$ of the operator $\mathcal A$ to these arcs satisfies
\begin{align*}
  \mathcal Bf(\zb\xi(\varphi,\vartheta))
  = \mathcal Af \left(Q\left(\tfrac{\vartheta}{2}, \tfrac{\pi}{2}, \tfrac{3\pi}{2}-\varphi\right), \tfrac{\vartheta}{2}\right)
  = \int_{0}^{\vartheta} f(\zb\eta(\varphi,\varrho)) \d\varrho.
\end{align*}
If $f $ is additionally {differentiable}, it can be recovered from $\mathcal Bf$ by
\begin{equation*}
f(\zb\xi(\varphi,\vartheta))
= \frac{\mathrm d}{\mathrm d\vartheta} \mathcal Bf(\zb\xi(\varphi,\vartheta)).
\end{equation*}

The following more general result for injectivity is due to \cite[Theorem 4.4.1]{Ami07}.
Its proof uses a similar idea combined with an extension by density.

\begin{Proposition}
	Let $S$ be an open subset of $\S^2$ and $A,B\subset S$ nonempty sets with $\overline{A \cup B} = \overline S$.
	If $f\in C(\S^2)$ and
	\begin{equation*}
	\int_{\gamma(\zb\xi,\zb\zeta)} f(\zb\eta) \d\zb\eta = 0
	\qquad \text{for all }\zb\xi\in A,\,\zb\zeta\in B,
	\end{equation*}
	then $f\equiv0$ on $S$.
\end{Proposition}

For $A = \{\zb e^3\}$ and $B=\S^2$, we have the arcs starting in the north pole.

\subsection{Recovery of local functions}
A subset $\Omega\subset\S^2$ is called convex if for any two points $\zb\xi,\zb\eta\in\Omega$ the geodesic arc $\gamma(\zb\xi,\zb\eta)$ is contained in $\Omega$.
We denote by $\partial\Omega$ the boundary of $\Omega$.

\begin{Theorem}
	Let $f\in C(\S^2)$ and $\Omega$ be a convex subset of $\S^2$ whose closure $\overline\Omega$ is strictly contained in a hemisphere, i.e., there exists a $\zb\zeta\in\S^2$ such that $\left<\zb\xi,\zb\zeta\right> >0$ for all $\zb\xi\in\overline\Omega$.
	If
	\begin{equation}
	\int_{\gamma(\zb\xi,\zb\eta)} f(\zb\eta) \d\zb\eta = 0
	\qquad\text{for all } \zb\xi,\zb\eta\in \partial \Omega,
	\label{eq:circle-arcs-vanish}
	\end{equation}
	then $f=0$ on $\Omega$.
	\label{thm:local}
\end{Theorem}
\begin{Proof}
	Without loss of generality, we assume that $\overline\Omega$ is strictly contained in the northern hemisphere, i.e., we have $\xi_3 >0$ for all $\zb\xi\in\overline\Omega$.
	We define the restriction of $f$ to $\overline\Omega$ by
	$$f_\Omega(\zb\xi)=
	\begin{cases}
	f(\zb\xi), &\zb\xi\in\overline\Omega\\
	0, &\zb\xi\in\S^2\setminus\overline\Omega
	\end{cases}$$
	Since $\gamma(\zb\xi,\zb\eta)\subset \overline\Omega$ for all $\zb\xi,\zb\eta\in \partial \Omega$,
	the function $f_\Omega$
	also satisfies \eqref{eq:circle-arcs-vanish}.

	For $\zb\xi\in\S^2$, denote with $\zb\xi^\perp = \{\zb\eta\in\S^2: \left<\zb\xi,\zb\eta\right>=0 \}$ the great circle perpendicular to $\zb\xi$.
	We show that the Funk--Radon transform
	\begin{equation}
	\mathcal Ff_\Omega (\zb\xi)
	= \int_{\zb\xi^\perp\cap\overline\Omega} f_\Omega(\zb\eta) \d \zb\eta
	+ \int_{\zb\xi^\perp\setminus\overline\Omega} f_\Omega(\zb\eta) \d \zb\eta
	\label{eq:arc-funk}
	\end{equation}
	vanishes everywhere.
	The second summand of \eqref{eq:arc-funk} vanishes because $f_\Omega$ is zero outside $\overline\Omega$ by definition.
	If $\zb\xi^\perp\cap\overline\Omega$ is not empty, there exist two points $\zb\eta^1,\zb\eta^2\in \partial\Omega$ such that $\gamma(\zb\eta^1,\zb\eta^2) = \zb\xi^\perp\cap\overline\Omega$, which shows that also the first summand of \eqref{eq:arc-funk} vanishes.
	Hence, $\mathcal Ff_\Omega=0$ on $\S^2$.
	Since the Funk--Radon transform $\mathcal F$ is injective for even functions, we see that $f_\Omega$ must be odd.
	Since $f_\Omega$ is supported strictly inside the northern hemisphere, so $f_\Omega$ must be the zero function.
	By the construction, we see that $f(\zb\xi)$ vanishes for all $\zb\xi\in\Omega$.
\end{Proof}

An analogue to \prettyref{thm:local} for $\Omega$ being the northern hemisphere and the arcs being half circles is shown in \cite{Rub16}.

\subsection{Arcs with fixed length}
\label{sec:fixed-length}

In the following, we consider circle arcs with fixed length $\psi$.
To this end, we define the restriction
\begin{equation*}
\mathcal A_\psi(Q) = \mathcal A(Q,\psi).
\end{equation*}

\begin{theorem}
	\label{thm:svd-fixed-length}
	Let $\psi\in(0,\pi)$ be fixed.
	The operator $\mathcal A_\psi\colon L^2(\S^2)\to L^2(\SO)$ has the singular value decomposition
	$$
	\left\{
	\left(Y_n^k, Z_{n,\psi}^k, \mu_n(\psi)\right):
	n\in\N,\,k\in\{-n,\dots,n\} \right\},
	$$
	with the singular values
	\begin{equation}
	\mu_{n}(\psi) = \sqrt{\sum_{j=-n}^n \frac{8\pi^2}{2n+1} \abs{\widetilde P_{n}^{j}(0)}^2 s_j(\psi)^2}
	\label{eq:mu}
	\end{equation}
	and the singular functions
	$$
	Z_{n,\psi}^k
	=  \frac{\mathcal A_\psi Y_n^k}{\mu_n(\psi)}
	= \frac{1}{\mu_n(\psi)} \sum_{j=-n}^{n} \widetilde P_n^j(0)\, s_j(\psi)\, D_n^{j,k}.
	$$
	In particular, $\mathcal A_\psi$ is injective.
\end{theorem}
\begin{Proof}
	Let $\psi\in(0,\pi)$ be fixed, $n\in\N$ and $k\in \{-n,\dots,n\}$.
	We have by \eqref{eq:AY}
	\begin{align*}
	&\left<\mathcal A_\psi Y_n^k, \mathcal A_\psi Y_{n'}^{k'}\right>_{L^2(\SO)}
	\\
	&\qquad= \sum_{j=-n}^n \sum_{j'=-n'}^{n'} \int_{\SO}
	D_n^{j,k}(Q)\, \overline{D_{n'}^{j',k'}(Q)}\, \widetilde P_{n}^{j}(0)\, \widetilde P_{n'}^{j'}(0)\, s_j(\psi)\, s_{j'}(\psi) \d Q
	\\&\qquad=
	\sum_{j=-n}^n \sum_{j'=-n'}^{n'} \frac{8\pi^2}{2n+1}\delta_{nn'}\,\delta_{kk'}\,\delta_{jj'}\, \widetilde P_{n}^{j}(0)\, \widetilde P_{n'}^{j'}(0)\, s_j(\psi)\, s_{j'}(\psi)
	\\&\qquad=
	\delta_{nn'}\delta_{kk'} \sum_{j=-n}^n \frac{8\pi^2}{2n+1} \abs{\widetilde P_{n}^{j}(0)}^2 s_j(\psi)^2 .
	\end{align*}
	For the injectivity, we check that the singular values $\mu_n(\psi)$
	do not vanish for each $n\in\N$.
	We have $\widetilde P_{n}^{j}(0) = 0$ if and only if $n-j$ is odd.
	Furthermore, the definition of $s_j$ in \eqref{eq:sj} shows that $s_0(\psi) =2\psi$ vanishes if and only if $\psi=0$ and $s_1(\psi) = 2\sin(\psi)$ vanishes if and only if $\psi$ is an integer multiple of $\pi$.
	Hence, the functions $\mathcal A Y_n^k$ are also orthogonal in the space $L^2(\SO)$.
\end{Proof}

\begin{Theorem}
	\label{thm:mu-lim}
	The singular values $\mu_n(\psi)$ of $\mathcal A_\psi$ satisfy for odd $n=2m-1$
	\begin{equation}
	\lim_{m\to\infty}
	\frac{4m-1}{4}\,\mu_{2m-1}(\psi)^2
	=
	\begin{cases}
	4\pi\psi, &\psi\in[0,\frac\pi2]\\
	4\pi^2-4\pi\psi, &\psi\in[\frac\pi2,\pi],
	\end{cases}
	\label{eq:mu-lim-odd}
	\end{equation}
	and for even $n=2m$
	\begin{equation}
	\lim_{m\to\infty}
	{\frac{4m+1}{4}}\,\mu_{2m}(\psi)^2
	=
	\begin{cases}
	4\pi\psi, &\psi\in[0,\frac\pi2]\\
	12\pi\psi-4\pi^2, &\psi\in[\frac\pi2,\pi].
	\end{cases}
	\label{eq:mu-lim-even}
	\end{equation}
\end{Theorem}
\begin{Proof}
We first show \eqref{eq:mu-lim-odd}.
Let $m\in\mathbb N$.
We have by \eqref{eq:mu}
\begin{align*}
\frac{4m-1}{4}\,\mu_{2m-1}(\psi)^2
&
= 16\pi^2 \sum_{j=1}^m \abs{\widetilde P_{2m-1}^{2j-1}(0)}^2 \frac{\sin^2((2j-1)\psi)}{(2j-1)^2}
.
\end{align*}
We denote by $\nu(\psi) = 4\pi\left(\frac\pi2-\abs{\psi-\frac\pi2}\right)$ the right-hand side of \eqref{eq:mu-lim-odd}.
The Fourier cosine series of $\nu$ reads by \cite[1.444]{GrRy07}
\begin{equation*}
16\sum_{k=1}^{\infty} \frac{\sin((2k-1)\psi)^2}{(2k-1)^2}
=16\sum_{k=1}^\infty \frac{1-\cos((2k-1)2\psi)}{2(2k-1)^2}
= \nu(\psi)
,\qquad \psi\in[0,\pi].
\end{equation*}
We have
\begin{align}
\norm{\frac{4m-1}{4}\mu_{2m-1}^2-\nu}_{C([0,\pi])}
&= \norm{16\sum_{j=1}^{\infty}\frac{\pi^2\abs{\widetilde P_{2m-1}^{2j-1}(0)}^2 - 1}{(2j-1)^2} \sin^2((2j-1)\psi) }_{C([0,\pi])}\nonumber
\\&
\le \sum_{j=1}^{\infty} \frac{16\,\abs{\pi^2\abs{\widetilde P_{2m-1}^{2j-1}(0)}^2 - 1}}{(2j-1)^2}.
\label{eq:mu-nu}
\end{align}
We show that \eqref{eq:mu-nu} goes to zero for $m\to\infty$, which then implies \eqref{eq:mu-lim-odd}.
By \eqref{eq:Pnj0-lim}, we see that $\pi^2\,\abs{\smash{\widetilde P_{2m-1}^{2j-1}(0)}}^2$ converges to $1$ for $m\to\infty$.
Using the singular values \eqref{eq:sv} together with their bound \eqref{eq:sv-bound-odd}, we obtain the following summable majorant of \eqref{eq:mu-nu}:
\begin{equation*}
\sum_{j=1}^{\infty} \frac{16\pi^2\,\abs{\widetilde P_{2m-1}^{2j-1}(0)}^2}{(2j-1)^2}
\le \frac{4m-1}{2\pi}\,\sigma_{2m-1}^2
\le 2\pi\,\frac{4m-1}{m} \left(\frac4{\sqrt3}+1\right)
.
\end{equation*}
Hence, the sum \eqref{eq:mu-nu} converges to 0 for $m\to\infty$ by the dominated convergence theorem of Lebesgue.

	In the second part, we show \eqref{eq:mu-lim-even} for the odd singular values.
	Let $m\in\mathbb N$.
	We have
	\begin{equation}
	\frac{4m+1}{4}\,\mu_{2m}(\psi)^2
	= 8\pi^2 \abs{\widetilde P_{2m}^{0}(0)}^2 \psi^2 + 4\pi^2 \sum_{k=1}^m \abs{\widetilde P_{2m}^{2k}(0)}^2 \frac{\sin^2(2k\psi)}{k^2}.
	\label{eq:mu2m}
	\end{equation}
	We examine both summands on the right side of \eqref{eq:mu2m}.
	The first summand converges due to \eqref{eq:Pnj0-lim}:
	\begin{equation*}
	\lim_{m\to\infty} 8\pi^2 \abs{\widetilde P_{2m}^{0}(0)}^2 \psi^2
	=8\psi^2.
	\end{equation*}
	We denote the second summand of \eqref{eq:mu2m} by
	\begin{equation*}
	\lambda_m(\psi) =
	4\pi^2 \sum_{k=1}^m \abs{\widetilde P_{2m}^{2k}(0)}^2 \frac{\sin^2(2k\psi)}{k^2}
	\end{equation*}
	and define $\lambda$ by the following Fourier cosine series, see \cite[1.443]{GrRy07},
	\begin{equation*}
	\lambda(\psi)=
	\sum_{k=1}^{\infty} \frac{\sin^2(2k\psi)}{k^2}
	=\sum_{k=1}^\infty \frac{1-\cos(4k\psi)}{2k^2}
	= \begin{cases}
	-2\psi^2+\pi\psi,& \psi\in[0,\frac\pi2)\\
	-2\psi^2+3\pi\psi-\pi^2,& \psi\in[\frac\pi2,\pi).
	\end{cases}
	\end{equation*}
	We have
	\begin{align*}
	\norm{\lambda_m - \lambda}_{C([0,\pi])}
	={}& \norm{\sum_{k=1}^{\infty}\frac{\pi^2\abs{\widetilde P_{2m}^{2k}(0)}^2 - 1}{k^2} \sin^2(2k\psi) }_{C([0,\pi])}
	\\
	\le{}& \sum_{k=1}^{\infty} \frac{\abs{\pi^2\abs{\widetilde P_{2m}^{2k}(0)}^2 - 1}}{(2j-1)^2}.
	\end{align*}
	As in the first part of the proof, we see with \eqref{eq:sv-bound-even} that the last sum goes to 0 for $m\to\infty$, which proves \eqref{eq:mu-lim-even}.
\end{Proof}
\begin{figure}

	\usetikzlibrary{intersections,calc}
	\tikzset{axis break gap/.initial=1mm}
	\begin{tikzpicture}
	\begin{axis}[
		name=bottom axis,
		legend cell align=left,
		xlabel= polynomial degree $n$ ,
		cycle list name=exotic,
		width=.99\textwidth,
		height=5.9cm,
		xmin = 0, xmax = 50,
		ymin = 0, ymax = 55,
		ytick={0,10,20,30,40,50},
		grid = major,
		axis x line*=bottom
		]
	\addplot+[thin,mark options={scale=1.0,solid}] table[x index=0,y index=5] {mu.dat};
	\addplot+[thin,mark options={scale=1.0,solid}] table[x index=0,y index=25] {mu.dat};
	\addplot+[thin,mark options={scale=1.0,solid}] table[x index=0,y index=45] {mu.dat};
	\addplot+[thin,mark options={scale=1.5,solid}] table[x index=0,y index=55] {mu.dat};
	\addplot+[thin,mark options={scale=1.0,solid}] table[x index=0,y index=98] {mu.dat};
	\end{axis}
	\begin{axis}[
		at=(bottom axis.north),
		name=top axis,
		anchor=south, yshift=\pgfkeysvalueof{/tikz/axis break gap},
		cycle list name=exotic,
		width=\textwidth,
		height=2.8cm,
		xmin = 0, xmax = 50,
		ymin = 145, ymax = 160,
		ytick={150,160},
		grid = major,
		legend entries = {$\psi=0.05\pi$,$\psi=0.25\pi$,$\psi=0.45\pi$,$\psi=0.55\pi$,$\psi=0.98\pi$},
		legend style={yshift=-3pt},
		axis x line*=top,
		xticklabel=\empty,
		after end axis/.code={
			\draw (rel axis cs:0,0) +(-2mm,-1mm) -- +(2mm,1mm)
			++(0pt,-\pgfkeysvalueof{/tikz/axis break gap})
			+(-2mm,-1mm) -- +(2mm,1mm)
			(rel axis cs:1,0) +(-2mm,-1mm) -- +(2mm,1mm)
			++(0pt,-\pgfkeysvalueof{/tikz/axis break gap})
			+(-2mm,-1mm) -- +(2mm,1mm);
		}]
	\addplot+[thin,mark options={scale=1.0,solid}] table[x index=0,y index=5] {mu.dat};
	\addplot+[thin,mark options={scale=1.0,solid}] table[x index=0,y index=25] {mu.dat};
	\addplot+[thin,mark options={scale=1.0,solid}] table[x index=0,y index=45] {mu.dat};
	\addplot+[thin,mark options={scale=1.5,solid}] table[x index=0,y index=55] {mu.dat};
	\addplot+[thin,mark options={scale=1.0,solid}] table[x index=0,y index=98] {mu.dat};
	\end{axis}
	\path (bottom axis.south west -| bottom axis.outer south west) -- node[anchor=south,rotate=90]{$(n+\frac12)\,\mu_n(\psi)^2$} (top axis.north west -| top axis.outer north west);
	\end{tikzpicture}

	\centering
	\begin{tikzpicture}
	\begin{axis}[xlabel=$\psi$, ylabel=$(n+\frac12)\,\mu_n(\psi)^2$,
	width=.8*\textwidth,height=0.56\textwidth,xmin=0, xmax=3.141592653589793,
	xtick={0,1.570796326794897,3.141592653589793},
	xticklabels={0,$\pi/2$,$\pi$},
	enlargelimits=false, legend pos=north west, legend style={cells={anchor=west}}]
	\addplot [mark=none, thick,smooth, samples=401, color=yellow!80!black,densely dashdotted] {4*pi*x^2};
	\addplot [mark=none, thick,smooth, samples=401, color=lime!80!black] {12*pi*sin((180*x)/pi)^2};
	\addplot [mark=none, thick,smooth, samples=401, color=lime!80!black,densely dashdotted] {(5*pi*(4*x^2 + 3*sin(180/pi*2*x)^2))/4};
	\addplot [mark=none, thick,smooth, samples=401, color=red] {4*pi^2*((21*sin((180*x)/pi)^2)/(8*pi) + (35*sin((540*x)/pi)^2)/(72*pi))};
	\addplot [mark=none, thick,smooth, samples=401, color=red,densely dashdotted] {(9*pi*(144*x^2 + 80*sin(180/pi*2*x)^2 + 35*sin(180/pi*4*x)^2))/256};
	\addplot [mark=none, thick,smooth, samples=401, color=blue] {	(11*pi*(6750*sin(180/pi*x)^2 + 875*sin(180/pi*3*x)^2 + 567*sin(180/pi*5*x)^2))/7200};
	\addplot [mark=none, thick,smooth, samples=401, color=blue,densely dashdotted] {	(13*pi*(1200*x^2 + 630*sin(180/pi*2*x)^2 + 189*sin(180/pi*4*x)^2 + 154*sin(180/pi*6*x)^2))/3072};
	\legend{$n=0$,$n=1$,$n=2$,$n=3$,$n=4$,$n=5$,$n=6$};
	\end{axis}
	\end{tikzpicture}

	\protect\caption[]{\label{fig:sv-fixed-length} The (normalized) singular values $(n+\frac12)\,\mu_n(\psi)^2$.
	Top: dependency on the degree $n$. Note the oscillation for $\psi>\frac\pi2$.
	Bottom: dependency on the arc-length $\psi$ (dashed lines correspond to even $n$).}
\end{figure}
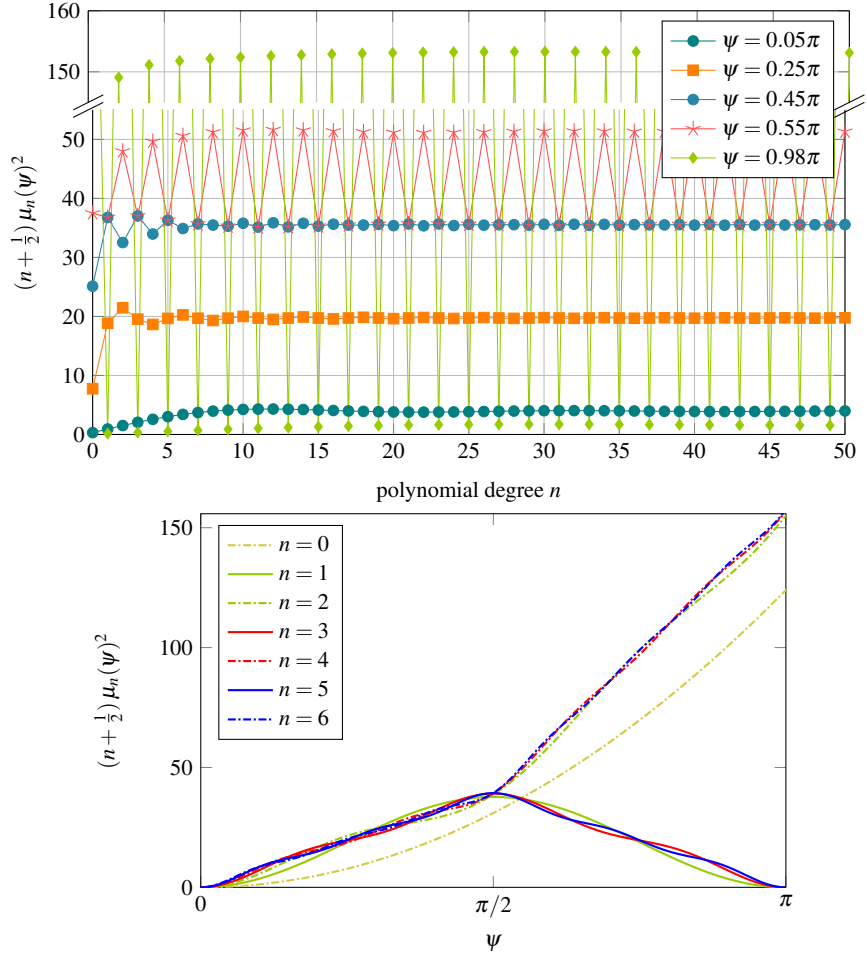

\begin{Remark}
Theorem \ref{thm:mu-lim} shows that the singular values $\mu_n(\psi)$ decay with the same asymptotic rate of $n^{-1/2}$ as $\sigma_n$ from Theorem \ref{thm:svd}.
For $\psi<\frac\pi2$, the singular values $(n+\frac12)\,\mu_n(\psi)^2$ for even and odd $n$ converge to the same limit.
However, for $\psi>\frac\pi2$, the singular values for even $n$ become larger than the ones for odd $n$.
This might be explained by the fact that for odd $n$, the spherical harmonics $Y_n^k$ are odd and integrating them along a circle arc with length $2\psi$, which is longer than a half-circle, yields some cancellation.
In the limiting case $\psi=\pi$, which is not covered by Theorem \ref{thm:svd-fixed-length}, $\mathcal A_{\pi}$ corresponds to the Funk--Radon transform, which is injective only for even functions and vanishes on odd functions.
This behavior is illustrated in Fig.\ \ref{fig:sv-fixed-length}.
\end{Remark}

\begin{Remark}
Since the rotation group $\SO$ is three-dimensional, the inversion of the arc transform $\mathcal A_\psi$ with fixed length is still overdetermined.

In the case $\psi=\frac\pi2$, we have the integrals along all half circles.
The injectivity of the arc transform for half circles was shown in \cite{Gro98}.
The restriction of the arc transform to all half circles that are subsets of either the upper or the lower hemisphere is still injective, see \cite{Rub16}.
This is because every function that is supported in the upper (lower) hemisphere can be uniquely reconstructed by its Funk--Radon transform, which then integrates only over the half circles in the upper (lower) hemisphere.
\end{Remark}

The singular value decomposition from \prettyref{thm:svd-fixed-length} allows us to reconstruct a function $f\in L^2(\S^2)$ given $g = \mathcal A_\psi f $.
\begin{Theorem}
\label{thm:inversion-fixed-length}
Let $f\in L^2(\S^2)$ and $g = \mathcal A_\psi f \in L^2(\SO)$.
Then $f$ can be reconstructed from the rotational Fourier coefficients $\hat g_n^{j,k}$ given in \eqref{eq:Fourier-coefficients-SO3} by
\begin{equation}
f = \sum_{n=0}^{\infty} \sum_{k=-n}^n
\frac{\sum_{j=-n}^{n} \widetilde P_n^j(0)\, s_j(\psi)\, \hat g_n^{j,k}}
{\sum_{j=-n}^{n} \widetilde P_n^j(0)^2\, {s_j(\psi)}^2}\, Y_n^k.
\label{eq:inversion-fixed-length}
\end{equation}
\end{Theorem}
\begin{Proof}
We have by Theorem \ref{thm:svd-fixed-length} for the spherical Fourier coefficients
\begin{align*}
\hat f_n^k
&= \frac{1}{\mu_n(\psi)} \left<g, Z_{n,\psi}^k\right>_{L^2(\SO)}
\\
&= \frac{1}{\mu_n(\psi)^2} \sum_{j=-n}^{n} \widetilde P_n^j(0)\, s_j(\psi) \left<g, D_n^{j,k} \right>_{L^2(\SO)}.
\end{align*}
The assertion follows by \eqref{eq:Fourier-coefficients-SO3} and \eqref{eq:mu}.
\end{Proof}

\begin{Remark}
	A big advantage of using the singular value decomposition for inversion is that it is straightforward to
	apply Tikhnov-type regularization or the mollifier method \cite{LoMa90}, which
	both correspond to a multiplication of the summands in the inversion formula
	\eqref{eq:inversion-fixed-length} with some filter coefficients $c_n$, cf.\ \cite{HiQu15}.
	We obtain
	\begin{equation}
	f_c = \sum_{n=0}^{\infty} \sum_{k=-n}^n c_n
	\frac{\sum_{j=-n}^{n} \widetilde P_n^j(0)\, s_j(\psi)\, \hat g_n^{j,k}}
	{\sum_{j=-n}^{n} \widetilde P_n^j(0)^2\, {s_j(\psi)}^2}\, Y_n^k.
	\label{eq:inversion-fixed-length-reg}
	\end{equation}
	Filter coefficients corresponding to Pinsker estimators are optimal for functions in certain Sobolev spaces, cf.\ \cite{FoGiHoHo16}.
	They were applied to the Funk--Radon transform in \cite{HiQu15}.
\end{Remark}

\section{Numerical tests}\label{sec:numeric}

We consider the arc transform $\mathcal A_\psi$ with fixed length $\psi\in(0,\pi)$ as in \prettyref{sec:fixed-length}.

\subsection{Forward algorithm}

For given $f\in C(\S^2)$, we want to compute the arc transform
$
\mathcal A_\psi f (Q_m)
$
at points $Q_m\in\SO$, $m=1,\dots,M$.
In order to derive an algorithm, we assume that $f\in \mathscr P_N(\S^2)$ is a polynomial.
We compute the spherical Fourier coefficients
$$
\hat f_n^k = \int_{\S^2} f(\zb\xi)\,\overline{Y_n^k(\zb\xi)}\d\zb\xi
,\qquad n=0,\dots,N,\; k=-n,\dots,n,
$$
with a quadrature rule on $\S^2$ that is exact for polynomials of degree $2N$.
The computation of the spherical Fourier coefficients $\hat f_n^k$ can be done with the adjoint NFSFT (Nonequispaced Fast Spherical Fourier Transform) algorithm \cite{KePo06} in $\mathcal O(N^2\,\log^2 N^2+M)$ steps.
Then, by \eqref{eq:AY},
\begin{equation}
\mathcal A_\psi f (Q_m)
= \sum_{n=0}^{N} \sum_{j,k=-n}^{n} \hat f_n^k\, \widetilde P_{n}^{j}(0)\, s_j(\psi)\, D_n^{j,k}(Q)
,\qquad m=1,\dots,M,
\label{eq:Af-discrete}
\end{equation}
is a discrete rotational Fourier transform of degree $N$,
which can be computed with the NFSOFT (Nonequispaced Fast $\SO$ Fourier Transform) algorithm \cite{poprvo07} in $\mathcal O(M+N^3\,\log^2 N)$ steps.
Implementations of both NFSFT and NFSOFT are contained in the NFFT library \cite{nfft33}.

A simple alternative for the computation of $\mathcal A_\psi f$ is the following quadrature with $K$ equidistant nodes
\begin{equation}
\mathcal A_\psi f(Q)
\approx \frac{2\psi}{K}\sum_{i=1}^{K}
f\left(Q^{-1}\left(\zb e_{\varrho_i}\right)\right)
,\qquad \rho_i = \frac{2i-1-K}{K}\,\psi.
\label{eq:Af-quad}
\end{equation}
Computing $\mathcal A_\psi f(Q_m)$ for $m=1,\dots,M$ with the quadrature rule \eqref{eq:Af-quad} requires $\mathcal O(KM)$ operations.
Hence, for a high number $M$ of evaluation nodes, the NFSOFT-based algroithm is faster than the quadrature
based on \eqref{eq:Af-quad}.

\subsection{Inversion}

\begin{figure}
	\captionsetup[subfigure]{justification=centering}
  \begin{subfigure}{.33\textwidth}
    \includegraphics[width=.99\textwidth]{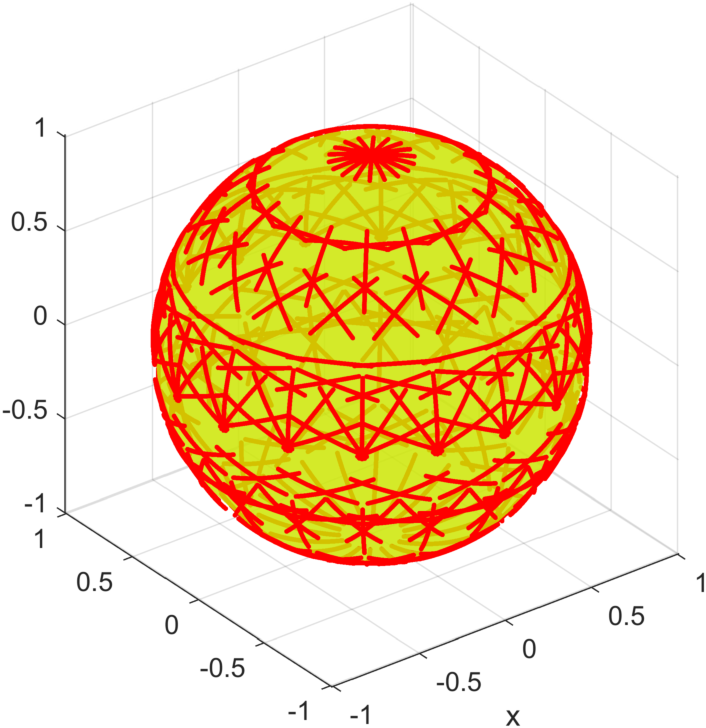}
    \caption{Gauss--Legendre\\ with 405 nodes}
    \label{fig:gl}
  \end{subfigure}\hfill%
  \begin{subfigure}{.33\textwidth}
    \includegraphics[width=.99\textwidth]{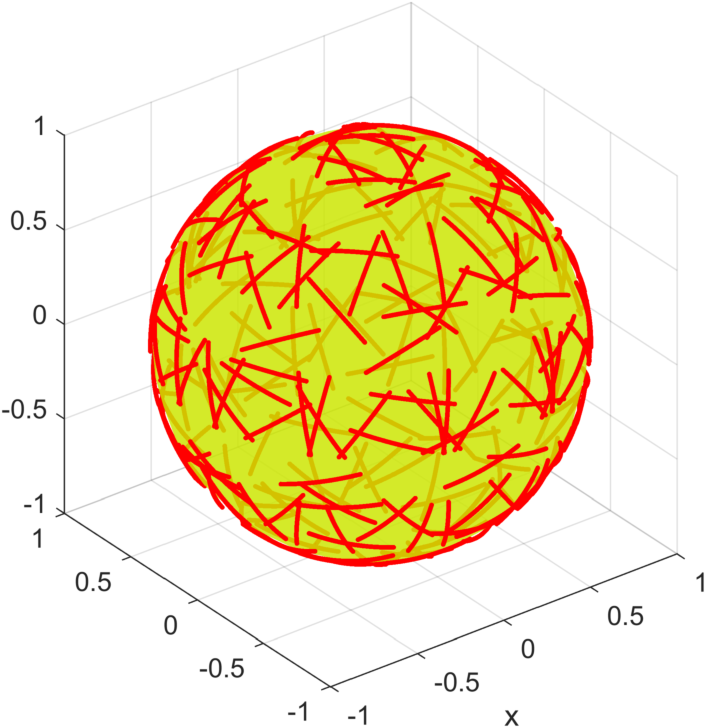}
    \caption{Tensor product $\S^1\times\S^2$\\ with 252 nodes}
    \label{fig:s1s2}
  \end{subfigure}\hfill%
  \begin{subfigure}{.33\textwidth}
    \includegraphics[width=.99\textwidth]{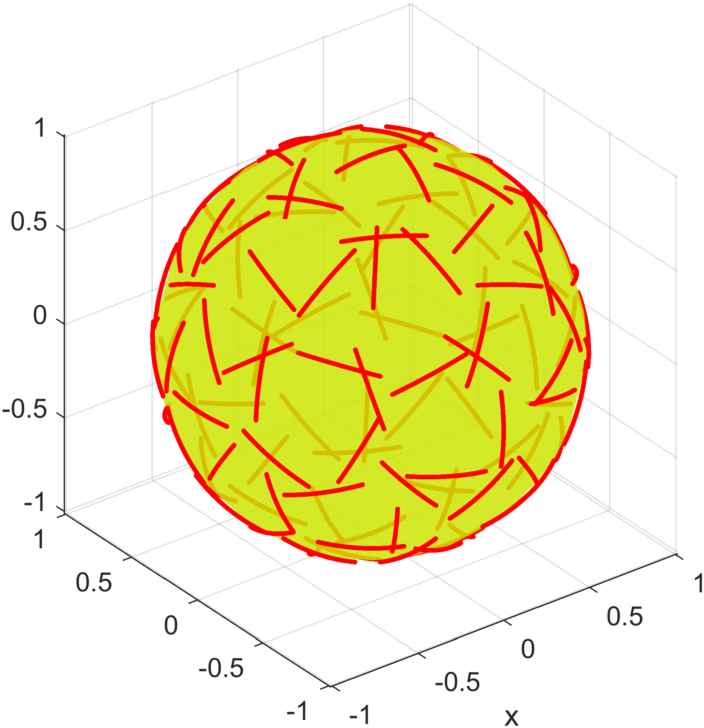}
    \caption{Gauss-type\\ with 240 nodes}
    \label{fig:mg}
  \end{subfigure}
  \caption{Circle arcs $\gamma(Q_m, 0.2)$ corresponding to quadrature nodes
    $Q_m\in\SO$, all quadrature formulas are exact for all rotational harmonics $D_n^{j,k}$ of
    degree $n\le8$.}
\label{fig:quadrature-arcs}
\end{figure}

We test the inversion from \prettyref{thm:inversion-fixed-length}. Let $g=\mathcal A_\psi f$. For the
computation of the rotational Fourier coefficients $\hat g_n^{j,k}$, $n=0,\dots,N$, $j,k=-n,\dots,n$, we use a
quadrature formula
\begin{equation}
  \hat g_n^{j,k}
  = \int_{\SO} g(Q)\, \overline{D_{n}^{j,k}(Q)} \d Q
   \approx  \sum_{m=1}^{M} w_m\, g(Q_m)\, \overline{D_{n}^{j,k}(Q_m)}
  \label{eq:quad-so}
\end{equation}
with nodes $Q_m\in\SO$ and weights $w_m>0$, $m=1,\dots,M$.
Again, we assume that $f\in\mathscr P_N(\S^2)$, which implies $\hat g_n^{j,k}=0$ for $n>N$.
Hence, \eqref{eq:quad-so} holds with equality if the quadrature integrates rotational harmonics up to degree $2N$ exactly.
There are different ways to obtain such exact quadrature formulas on $\SO$.
In a tensor product approach, we use Gauss--Legendre quadrature in $\cos\beta$ and a trapezoidal rule in both $\alpha$ and $\gamma$.
We can also write $\SO \sim \S^1\times\S^2$ and pair a trapezoidal rule on $\S^1$ in $\alpha$ with a quadrature on $\S^2$ with azimuth $\gamma$ and polar angle $\beta$, see \cite{GrPo09}.
Furthermore, Gauss-type quadratures on $\SO$ that are exact up to machine precision were computed in \cite{DissGr13}.
In \prettyref{fig:quadrature-arcs}, one can see the circle arcs corresponding to different quadrature rules on $\SO$, namely Gauss--Legendre nodes (Fig.\ \ref{fig:gl}), the tensor product of $\S^1\times\S^2$ (Fig.\ \ref{fig:s1s2}) and a Gauss-type quadrature on $\SO$ (Fig.\ \ref{fig:mg}).
We used Gauss-type quadratures on both $\S^2$ and $\SO$ from \cite{Gr_points}.
Note that because of the symmetry of the Gauss-type quadrature on $\SO$ we used in Fig\ \ref{fig:mg}, every arc corresponds to two quadrature nodes on $\SO$.

\begin{figure}
	\begin{subfigure}{.33\textwidth}
		\includegraphics[width=\textwidth]{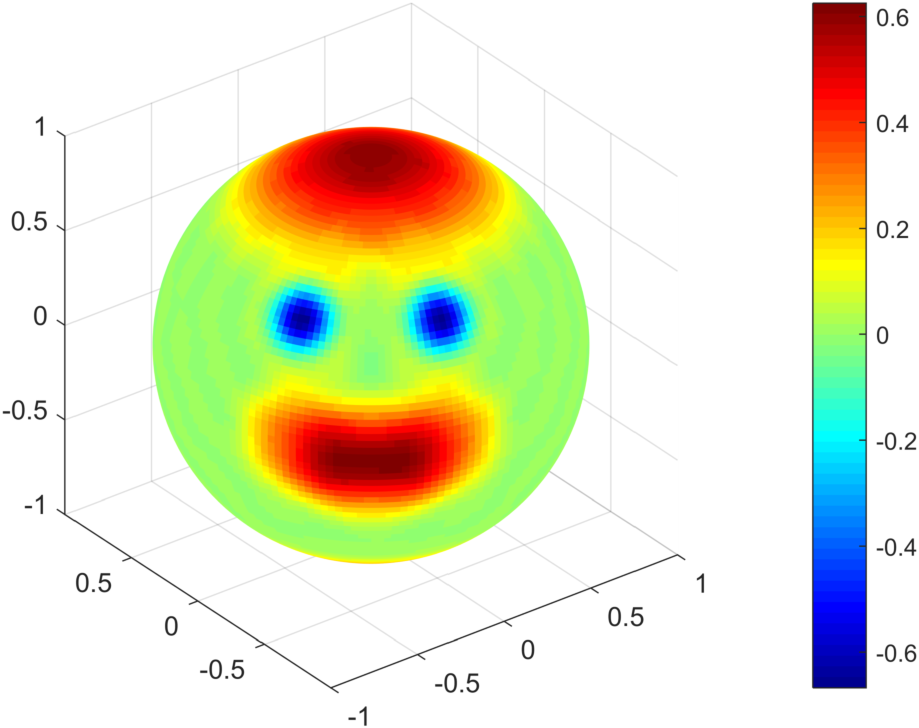}
		\caption{Exact data}
	\end{subfigure}
	\begin{subfigure}{.33\textwidth}
	\includegraphics[width=\textwidth]{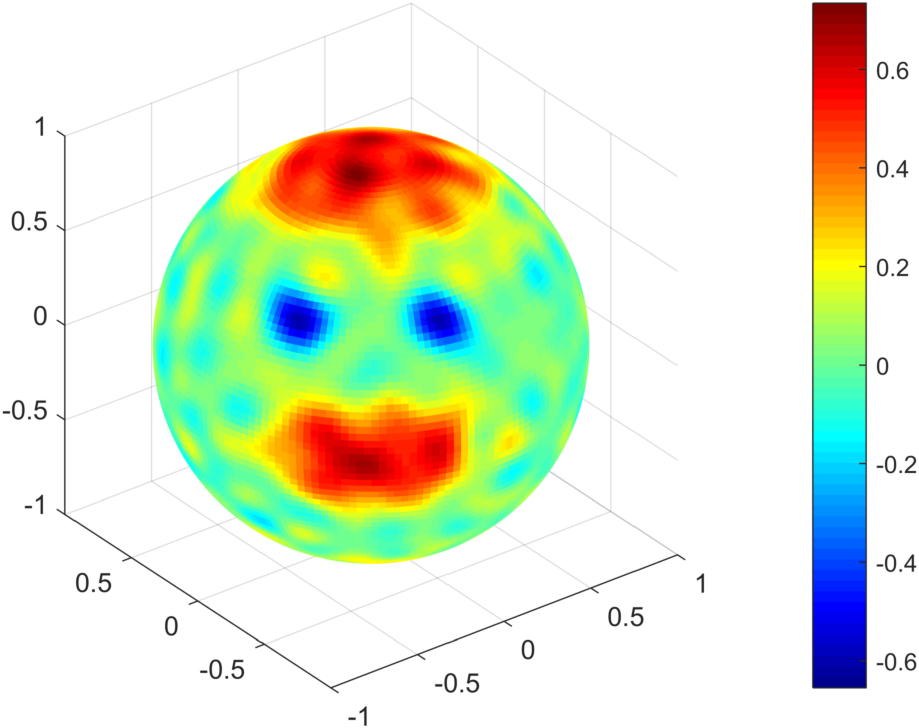}
	\caption{Noisy data}
	\end{subfigure}
	\begin{subfigure}{.33\textwidth}
		\includegraphics[width=\textwidth]{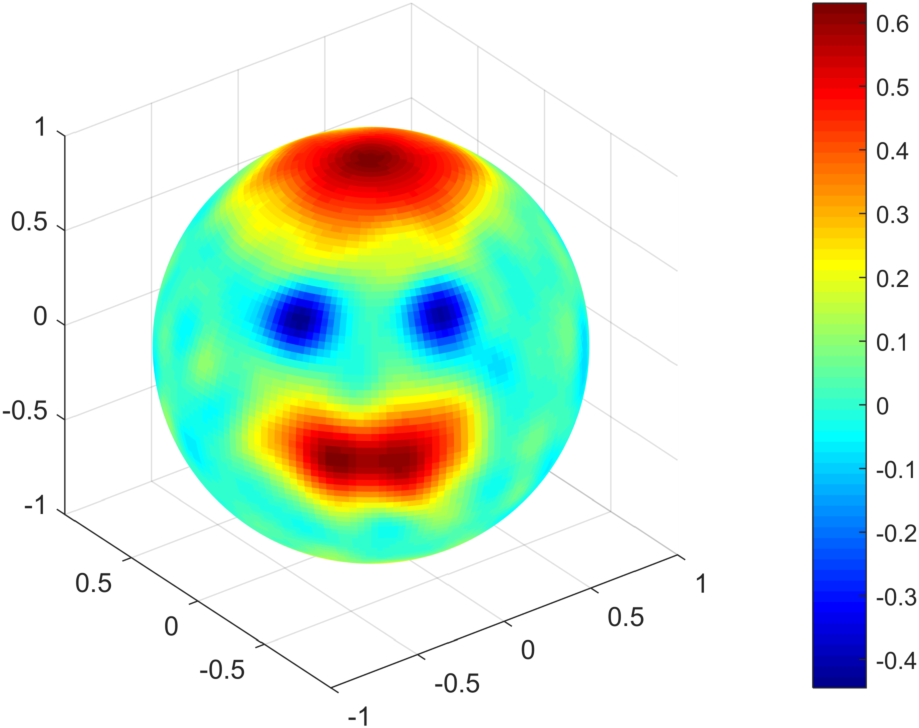}
		\caption{Noisy data \& regularization}
	\end{subfigure}
	\caption{Reconstruction of a spherical test function $f$ for degree $N=22$, $\psi=0.7$ and a tensor product $\S^1\times\S^2$ quadrature ($M=30240
		$).
}
	\label{fig:inversion}
\end{figure}

The reconstruction formula \eqref{eq:inversion-fixed-length} becomes the discrete rotational Fourier transform
\begin{equation*}
f = \sum_{n=0}^{N} \sum_{k=-n}^n
\frac{\sum_{j=-n}^{n} \widetilde P_n^j(0)\, s_j(\psi)\, \hat g_n^{j,k}}
{\sum_{j=-n}^{n} \widetilde P_n^j(0)^2\, {s_j(\psi)}^2}\, Y_n^k.
\end{equation*}
In Fig.\ \ref{fig:inversion}, we compare the reconstruction results,
where we use an artificial test function, the parameter $N=22$ and the tensor product of a trapezoidal rule on $\S^1$ with a Gauss-type quadrature on $\S^2$ from \cite{Gr_points}.
The resulting $\SO$ quadrature uses $M=30240$ nodes and is exact for degree 44.
We first perform the inversion without any noise in the data.
The reconstruction has an RMSE (root mean square error) of $0.0338$.
Then we add Gaussian white noise with a standard deviation of $0.2$ to the data $\mathcal A_\psi f(Q_m)$ and achieve an RMSE of $0.2272$.
Even though we did not perform any regularization, the reconstruction from
noisy data still looks considerably well.
This might be explained by the fact that the inverse arc transform with fixed
  opening angle and full $\SO$ data is still an overdetermined problem.
Applying the regularization \eqref{eq:inversion-fixed-length-reg} truncated to degree $n\le N$ with filter coefficients from \cite{HiQu15} yields a smaller RMSE of $0.1393$.

\begin{acknowledgement}
	The authors thank Volker Michel for pointing out the problem of spherical surface wave tomography at the Mecklenburg Workshop on Approximation Methods and Data Analysis 2016 and for fruitful conversations later on.
	Furthermore, we thank the anonymous reviewer for providing helpful comments and suggestions to improve this article.
\end{acknowledgement}

\end{document}